\documentclass[12pt,reqno]{amsart}%
\usepackage{amsfonts}
\usepackage{amsmath}
\usepackage{stmaryrd}
\usepackage{amssymb}
\usepackage[title]{appendix}
\usepackage{graphicx}%
\providecommand{\U}[1]{\protect\rule{.1in}{.1in}}
\newtheorem{theorem}{Theorem}[section]
\theoremstyle{plain}

\newtheorem{lemma}{Lemma}[section]

\newtheorem{proposition}{Proposition}
\newtheorem{remark}{Remark}

\numberwithin{equation}{section}
\begin{document}
\title[Sharp weighted bounds for reverse fractional integral operators]{Sharp weighted bounds for reverse fractional integral operators}
\author{Liang Huang}
\address{School of Science\\ Xi'an University of Posts and Telecommunications\\Xi'an, 710121, China\\ E-mail: huangliang10@163.com}
\author{Hanli Tang}
\address{Laboratory of Mathematics and Complex Systems (Ministry of Education), School of Mathematical Sciences\\Beijing Normal University\\Beijing, 100875, China\\ E-mail: hltang@bnu.edu.cn}
\thanks{This project is supported by National Key Research and Development Program of China(Grant No.2020YFA0712900) and National Science Foundation of China (Grant No. 12501127)
\\Corresponding author: Hanli Tang at hltang@bnu.edu.cn}


\begin{abstract}
In this paper, we prove the one-weight reverse fractional integral operators
\[
\|I_s f\|_{L^q(w^q)}\geq C_w\|f\|_{L^p(w^p)}
\]
holds for all nonnegative measurable functions $f$, with $\frac ns<p<1$ and
$\frac1q=\frac1p-\frac sn$, if and only if
$w^{p/(1-p)}\in A_{\frac{p(s-n)}{n(1-p)}}$.
Furthermore, we show the lower bound
\[
C_w\gtrsim_{n,s,p}[w]_{p,q}^{\min\left\{\frac{p(n-s)}{n(1-p)},\frac{p(n-s)}{sp-n}\right\}},
\]
and the exponent is sharp.

We also investigate when the reverse weak-type estimate
\[
\|I_s f\|_{L^{\frac{n}{n-s},\infty}(w^{\frac{n}{n-s}})}\geq \frac{n}{s}\,v_n^{(n-s)/n}[w^{\frac{n}{n-s}}]^{(n-s)/n}_{A_1}\|f\|_{L^1(w)}
\]
holds if and only if $w^{\frac{n}{n-s}}\in A_1$.
\end{abstract}

\maketitle

\section{introduction}
The classical fractional integral operators (or Riesz potentials) \(I_{\alpha}\) are defined by
\[
I_{\alpha}(f)(x)=\gamma_\alpha\int_{\mathbb{R}^n}\frac{f(x-y)}{|y|^{n-\alpha}}dy,
\]
where \(0<\alpha<n\) and \(\gamma_\alpha=2^{-\alpha}\pi^{-\frac{n}{2}}\frac{\Gamma(\frac{n-\alpha}{2})}{\Gamma(\frac{\alpha}{2})}\). This family of fractional potential operators was first systematically investigated by Riesz~\cite{R}. The \((L^p,L^q)\)-boundedness of Riesz potentials was proved by Hardy and Littlewood~\cite{HL} in the one-dimensional case, and by Sobolev~\cite{S} in higher dimensions. The endpoint weak-type \((L^1,L^{\frac{n}{n-\alpha},\infty})\) boundedness was subsequently obtained by Zygmund~\cite{Zygmund1956}. The classical result is the following.
\vskip0.2cm
\textbf{Theorem A. }\textit{Let} \(0<\alpha<n\) \textit{and let}
 \(p,q\) \textit{satisfy} \(1\leq p<q<\infty\) \textit{and} \(\frac{1}{p}-\frac{1}{q}=\frac{\alpha}{n}\). \textit{Then for} \(p>1\),
\[
\|I_{\alpha}(f)\|_{L^q(\mathbb{R}^n)}\leq{C(n,p,\alpha)}\|f\|_{L^{p}(\mathbb{R}^n)}.
\]
 \textit{When }$p=1$,
\[
\|I_\alpha(f)\|_{L^{\frac{n}{n-\alpha},\infty}(\mathbb{R}^n)}=\sup_{\lambda>0}\lambda|\{x\in{\mathbb{R}^n}:|I_{\alpha}f|>\lambda\}|^{\frac{n-\alpha}{n}}\leq{C(n,\alpha)\|f\|_{L^{1}(\mathbb{R}^n)}}.
\]
\vskip0.2cm
Lieb~\cite{L} calculated the sharp constant in the conformally symmetric case \(p=\frac{2n}{n+\alpha}\), \(q=\frac{2n}{n-\alpha}\) (see also \cite{FL}). More recently, Tang and Wang~\cite{TangWang2026} determined sharp weak-type constants for classical fractional integral operators.

The study of reverse fractional integral operators was initiated by Dou and Zhu~\cite{DouZhu2015}, who proved the existence of extremal functions and computed the best constant for $p=q=\frac{2n}{n+s}$. Chen, Liu, Lu and Tao~\cite{ChenLiuLuTao2018} established the reverse weighted Hardy-Littlewood-Sobolev inequality, namely the reverse Stein-Weiss inequality,
\[
\int_{\mathbb R^n}\int_{\mathbb R^n}
|x|^\alpha\,|x-y|^{s-n}\,f(x)\,g(y)\,|y|^\beta\,dx\,dy
\geq C_{n,\alpha,\beta,p,q'}\,\|f\|_{L^{q'}}\,\|g\|_{L^p}
\]
for any nonnegative functions $f\in L^{q'}(\mathbb R^n)$,
$g\in L^p(\mathbb R^n)$, and $p,q'\in(0,1)$, $\alpha,\beta,s-n>0$ such
that
\[
\frac1p+\frac1{q'}-\frac{\alpha+\beta+s-n}{n}=2.
\]
Kassymov, Ruzhansky and Suragan~\cite{KassymovRuzhanskySuragan2022} subsequently extended the reverse weighted and unweighted Hardy-Littlewood-Sobolev inequality to homogeneous groups. For the endpoint case, the present authors jointly with Zhang~\cite{HuangTangZhang2025} established the sharp reverse weak estimates and the extremal functions for the reverse fractional integral operators \(I_s\). For related reverse Hardy-Littlewood-Sobolev inequality, we refer the reader to \cite{CarrilloEtAl2019,NgoNguyen2017}.

In the weighted setting, Stein and Weiss \cite{SteinWeiss1958} established
the following weighted Hardy-Littlewood-Sobolev inequality (also called the Stein--Weiss inequality):
\begin{equation}\label{eq:stein-weiss}
\int_{\mathbb R^n}\int_{\mathbb R^n}
|x|^{-\alpha}|x-y|^{-\lambda}f(x)g(y)|y|^{-\beta}\,dx\,dy
\leq C_{n,\alpha,\beta,p,q}\,\|f\|_{L^p(\mathbb R^n)}\|g\|_{L^q(\mathbb R^n)},
\tag{1.3}
\end{equation}
where $1<p,q<\infty$, $\alpha,\beta,\lambda$ satisfy
\[
\frac1p+\frac1q+\frac{\alpha+\beta+\lambda}{n}=2,
\qquad
\frac1p+\frac1q\geq1,
\]
\[
\alpha+\beta\geq0,\qquad
\alpha<\frac{n}{p'},\qquad
\beta<\frac{n}{q'},\qquad
0<\lambda<n.
\]

The $A_p$ weight theory began with Muckenhoupt's theorem~\cite{Muckenhoupt1972Maximal}: the Hardy--Littlewood maximal operator \(M\) is bounded on \(L^p(w)\) for \(1<p<\infty\) if and only if the weight \(w\) satisfies the \(A_p\) condition
\begin{equation}\label{eq:Ap}
[w]_{A_p} := \sup_Q \left(\frac{1}{|Q|}\int_Q w \right)
\left(\frac{1}{|Q|}\int_Q w^{-1/(p-1)} \right)^{p-1} < \infty.
\end{equation}
Hunt, Muckenhoupt, and Wheeden~\cite{HuntMuckenhouptWheeden1973} proved that the
$A_p$ condition also characterizes the boundedness of the Hilbert transform
\[
Hf(x)=\mathrm{p.v.}\,\frac{1}{\pi}\int\frac{f(y)}{x-y}\,dy
\]
on $L^p(w)$. Coifman and Fefferman~\cite{CoifmanFefferman1974} developed the general theory for Calder\'on--Zygmund operators and introduced the class \(A_\infty = \bigcup_{p>1} A_p\).  Muckenhoupt and
Wheeden~\cite{MuckenhouptWheeden1974} characterized the boundedness of the fractional integral operators as follows:
\begin{equation}\label{eq:mw-strong}
\left(\int_{\mathbb R^n}|I_\alpha f(x)\,w(x)|^q\,dx\right)^{1/q}
\leq C\left(\int_{\mathbb R^n}|f(x)\,w(x)|^p\,dx\right)^{1/p}
\quad\text{for all }f\geq0,
\end{equation}
where $1/q=1/p-\alpha/n$, holds if and only if \([w]_{{p,q}}<\infty\), where the weight characteristic is defined as
\[
[w]_{{p,q}} =\sup_Q \left(\frac{1}{|Q|}\int_Q w^q dx\right)\left(\frac{1}{|Q|}\int_Q w^{-p'} dx\right)^{q/p'} <\infty.
\]
The $[w]_{{p,q}}$ condition was first introduced in~\cite{MuckenhouptWheeden1974}. Sawyer~\cite{Sawyer1982} first gave the two-weight
characterization for the maximal operator, and then set up the two-weight
bounds for fractional and Poisson integrals~\cite{Sawyer1988}. Sawyer and
Wheeden~\cite{SawyerWheeden1992} extended the two-weight characterization of
fractional integrals to Euclidean and homogeneous spaces. P\'erez established
two-weight norm inequalities for Riesz potentials~\cite{Perez1990}, and subsequently for potential operators and
fractional maximal operators~\cite{Perez1994}.

Once the weighted boundedness of the operator is obtained, it is natural to
consider how the operator bound depends on the weight. Buckley~\cite{Buckley1993} proved that for the maximal operator,
\[
\|M\|_{L^p(w)\to L^p(w)}\le c\,[w]_{A_p}^{1/(p-1)},
\]
and that the exponent \(1/(p-1)\) is sharp. For singular integral operators, the central open problem was so-called \(A_2\) conjecture, which asserts that for every Calder\'{o}n--Zygmund operator \(T\),
\[
\|T\|_{L^2(w)\to L^2(w)}\le C_T\,[w]_{A_2}.
\]
Petermichl~\cite{Petermichl2007,Petermichl2008} showed that
\[
\|T\|_{L^p(w)\to L^p(w)}\le c\,[w]_{A_p}^{\max\{1,1/(p-1)\}},
\]
where \(T\) is either the Hilbert transform or one of the Riesz transforms
\[
R_j f(x)=c_n\,\mathrm{p.v.}\int_{\mathbb R^n}
\frac{x_j-y_j}{|x-y|^{n+1}}f(y)\,dy.
\]
The \(A_2\) conjecture was resolved for general Calder\'on--Zygmund operators by Hyt\"onen~\cite{Hytonen2012}. This landmark proof was subsequently simplified by Lerner \cite{Lerner2013}, who gave an elegant argument via sparse pointwise domination. The sharp extrapolation theorem of Dragi\v{c}evi\'c, Grafakos, Pereyra, and Petermichl \cite{DGPP2005} extended the sharp \(p=2\) result to all \(1<p<\infty\), yielding
\[
\|T\|_{L^p(w) \to L^p(w)} \lesssim [w]_{A_p}^{\max(1, 1/(p-1))}.
\]
Furthermore, the sharp quantitative reverse H\"older theorem of Hyt\"onen, P\'erez, and Rela~\cite{HytonenPerezRela2012} and the sharp Fujii--Wilson \(A_\infty\) bounds of Hyt\"onen and P\'erez~\cite{HytonenPerez2013} have also been established. Recently, Lacey, Moen, P\'erez, and Torres~\cite{LaceyMoenPerezTorres2010} obtained the corresponding sharp bounds for weighted fractional integral operators \(I_\alpha\):
\vskip0.2cm
\textbf{Theorem B. }
Let \(1 < p < n/\alpha\) and \(q\) be defined by the equation \(1/q = 1/p - \alpha/n\), and \([w]_{{p,q}}<\infty\). Then,
\begin{equation*}
\|I_\alpha\|_{L^p(w^p)\to L^q(w^q)} \le c\,[w]_{{p,q}}^{\left(1-\frac{\alpha}{n}\right)\max\left\{1,\frac{p'}{q}\right\}}.
\end{equation*}
Furthermore this estimate is sharp.
\vskip0.2cm
Further one-weight and two-weight results for fractional integral operators were obtained by Cruz-Uribe and Moen~\cite{CruzUribeMoen2013}, who extend Theorem B to
\[
\|I_\alpha\|_{L^p(w^p)\to L^q(w^q)} \lesssim [w]_{A_p}^{1/p} \left( [w]_{A_\infty}^{1/p'} + [w^{1-p'}]_{A_\infty}^{1/p} \right).
\]

These existing results largely focus on power-weight settings, leaving open the general Muckenhoupt type one-weight characterization and sharp quantitative bounds for the best constant in reverse fractional integral operators. In the present paper we address two main questions for the reverse fractional integral operators \(I_s\) with $s>n$:

\begin{enumerate}
    \item What weight condition characterizes the boundedness of the reverse
    fractional integral operator, and what is the sharp dependence of the best
    constant on that weight?
    \item What is the corresponding weight condition for endpoint boundedness?
\end{enumerate}
The reverse fractional integral operators \(I_{s}\) are defined by
\[
I_{s}(f)(x)=\int_{\mathbb{R}^n}|x-y|^{n-s}f(y)dy,
\]
where \(s>n\). Throughout the strong-type discussion, set $p'=p/(p-1)<0$ and
$q'=q/(q-1)\in(0,1)$. Define the scale-invariant weight characteristic
\begin{equation*}
[w]_{p,q}:=\sup_B\left(\frac{1}{|B|}\int_B w^{-p'}dx\right)^{-\frac{1}{p'}}
\left(\frac{1}{|B|}\int_B w^{q}dx\right)^{-\frac{1}{q}}.
\end{equation*}
Here the supremum runs over all Euclidean balls $B\subset\mathbb{R}^n$. By direct computation,
\begin{equation}\label{eq:Arelation}
[w^{-p'}]_{A_{1+p'/q}}=[w]_{p,q}^{-p'},
\end{equation}
so that $[w]_{p,q}<\infty$ is equivalent to $w^{-p'}\in A_{1+p'/q}$.

Our first main theorem provides the qualitative one-weight characterization.
\begin{theorem}\label{thm:characterization}
Let $n\ge1$, $s>n$, $n/s<p<1$, and let $q<0$ satisfy $1/q=1/p-s/n$. Let $w>0$ almost everywhere. Then there exists a constant $C>0$ such that
\begin{equation}\label{eq:strong}
\|I_s f\|_{L^q(w^q)}\ge C\|f\|_{L^p(w^p)}
\end{equation}
for every nonnegative measurable $f$ if and only if $[w]_{p,q}<\infty$.
\end{theorem}
\begin{remark}
Let $w(x)=|x|^\gamma$. Then \eqref{eq:strong} holds if and only if
\begin{equation}\label{eq:power-range}
n-\frac np<\gamma<s-\frac np.
\end{equation}
In fact, by Theorem \ref{thm:characterization}, the required condition is \(|x|^{-p'\gamma}\in A_{1+p'/q}\), then we have
\[
-n<-p'\gamma<n(1+p'/q-1)=\frac{np'}{q}.
\]
Thus $n-\frac np=n/p'<\gamma<-n/q=s-\frac np$. The range \eqref{eq:power-range} agrees with the one-weight specialization of the reverse Stein-Weiss parameter range in~\cite{ChenLiuLuTao2018}. This gives an independent consistency check on the Muckenhoupt characterization.
\end{remark}
Next, we obtain the sharp weighted bounds for the reverse fractional integral operators:
\begin{theorem}\label{thm:mixed}
Let $n\ge1$, $s>n$, $n/s<p<1$, and let $q<0$ satisfy $1/q=1/p-s/n$ and \([w]_{p,q}<\infty\). Then,
\begin{equation*}
\|I_s\|_{L^p(w^p)\to L^q(w^q)} \ge c_{n,s,p}\,
[w]_{p,q}^{\min\left\{
\frac{p(n-s)}{n(1-p)},
\frac{p(n-s)}{sp-n}
\right\}}.
\end{equation*}
Furthermore this estimate is sharp.
\end{theorem}
\begin{remark}
Through the proof of Theorem \ref{thm:mixed}, we have the following better bounds,
\begin{equation}\label{eq:mixed-main}
C_*(w)\ge c_{n,s,p}\,
[w]_{p,q}^{-1}
\min\left\{[w^{-p'}]_{A_\infty}^{1/q},[w^q]_{A_\infty}^{1/p'}\right\}.
\end{equation}
Combining (\ref{eq:mixed-main}), (\ref{eq:Arelation}), the inequalities
$$[w^{-p'}]_{A_\infty}\lesssim_{n,p,q}[w^{-p'}]_{A_{1+p'/q}}\mbox{~and~} [w^q]_{A_{1+q/p'}}=[w^{-p'}]_{A_{1+p'/q}}^{q/p'},$$
we obtain the desired lower bound in Theorem \ref{thm:mixed}.
\end{remark}

The critical point is \(p_c=\frac{2n}{n+s}\), where the two symmetric terms coincide, giving \(p'=q\). Thus
\begin{equation*}
\Gamma_*=\begin{cases}
\displaystyle \frac{p(s-n)}{sp-n},&\displaystyle \frac ns<p\le \frac{2n}{n+s},\\[2mm]
\displaystyle \frac{p(s-n)}{n(1-p)},&\displaystyle \frac{2n}{n+s}\le p<1.
\end{cases}
\end{equation*}
At the conformally symmetric reverse Hardy-Littlewood-Sobolev exponent $p=2n/(n+s)$, the sharp exponent is
\[
C_*(w)\gtrsim [w]_{p,q}^{-2}.
\]

At the endpoint \(p=1\), we obtain a lower bound estimate of the reverse fractional integral operators.
\begin{theorem}\label{thm:endpoint-main}
Let $n\ge1$, $s>n$. The inequality
\[\|I_sf\|_{L^{\frac{n}{n-s},\infty}(w^{\frac{n}{n-s}})}\ge \frac ns v_n^{(n-s)/n}[w^{\frac{n}{n-s}}]_{A_1}^{(n-s)/n}\|f\|_{L^{1}(w)}\]
holds for all $f\ge0$ with a positive constant if and only if $w^{\frac{n}{n-s}}\in A_1$.
\end{theorem}

For two quantities $X,Y$ depending on parameters $p,s,n,\dots$,
the notation $X\lesssim Y$ (resp.\ $X\gtrsim Y$) means $X\leq C Y$ (resp.\ $X\geq c Y$)
for constants $c,C>0$ that depend only on the indicated subscripts;
$X\approx Y$ means both $X\lesssim Y$ and $X\gtrsim Y$.

The paper is organised as follows. Section~2 collects facts on $A_p$ and
$A_\infty$ weights, the bilinear reduction, and some lemmas.
Section~3 proves necessity of the Muckenhoupt condition. In Section~4,
a negative-index duality argument reduces the lower bound to a bilinear case, and we prove Theorems~\ref{thm:characterization} and~\ref{thm:mixed}. We use power-weight to show that the exponent is sharp. Section~5 treats the
endpoint $p=1$ and proves Theorem~\ref{thm:endpoint-main} by the layer-cake identity, the $A_1$ condition, and Tonelli's theorem,
while necessity follows by testing on approximate point masses and applying the  Lebesgue
differentiation.

\section{Preliminaries}\label{sec:Ap-background}
In this section, we recall facts on $A_r$ and $A_\infty$ weights and establish the bilinear and geometric lemmas used below.

In the main theorems, $w$ is positive and finite almost everywhere; local
integrability of $w$ itself is not required. We use $w(E)=\int_E w(x)\,dx$
when this integral is defined. For $1<r<\infty$, a positive, locally
integrable function $u$ belongs to $A_r$ if
\begin{equation*}
[w]_{A_r}
:=\sup_B\left(\frac{1}{|B|}\int_B w\right)
\left(\frac{1}{|B|}\int_B w^{-1/(r-1)}\right)^{r-1}<\infty.
\end{equation*}
Ball and cube definitions yield the same $A_r$ class, with comparison constants depending on the dimension and the exponent, see \cite{Muckenhoupt1972Maximal}.

The dual weight is
\[
w^*=w^{-1/(r-1)},
\]
and
\begin{equation}\label{eq:dualAp}
[w^*]_{A_{r'}}=[w]_{A_r}^{1/(r-1)}.
\end{equation}
Introducing the scale-invariant characteristic
\[
[w]_{p,q}=\sup_B\left(\frac{1}{|B|}\int_B w^{-p'}dx\right)^{-\frac{1}{p'}}
\left(\frac{1}{|B|}\int_B w^{q}dx\right)^{-\frac{1}{q}}.
\]
Equivalently, $[w]_{p,q}<\infty$ means
\[
w^{-p'}=w^{p/(1-p)}\in A_{1+p'/q}
=A_{\frac{p(s-n)}{n(1-p)}}.
\]

It is well known that $A_\infty=\bigcup_{r>1}A_r$. Equivalently, we use the Fujii-Wilson characteristic
\begin{equation}\label{eq:Ainfdef}
[w]_{A_\infty}
:=\sup_B \frac{1}{w(B)}\int_B M(w\chi_B)(x)\,dx,
\end{equation}
where $M$ is the uncentered Hardy--Littlewood maximal operator. Equivalent $A_\infty$ characteristics are discussed in, for example, \cite{Wilson2008,HytonenPerez2013}. The crucial quantitative fact is the sharp reverse-H\"older
theorem of Hyt\"onen--P\'erez--Rela~\cite{HytonenPerezRela2012}. In the
Euclidean setting it implies that there are dimensional constants
$c_n,C_n>0$ such that, for every $w\in A_\infty$, every ball $B$, and every
measurable $E\subset B$,
\begin{equation}\label{eq:Ainf-decay}
\frac{w(E)}{w(B)}\leq C_n\left(\frac{|E|}{|B|}\right)^{\delta_w},
\qquad
\delta_w\geq\frac{c_n}{[w]_{A_\infty}}.
\end{equation}

For $q<0$ and a positive measurable weight $w$, define
\begin{equation*}
\|f\|_{L^q(w)}
:=\left(\int_{\mathbb{R}^n} f(x)^q w(x)\,dx\right)^{1/q},
\end{equation*}
with the standard extended-value convention. The $L^p$ norm for $p<1$ is not a norm for a vector space. For $q<0$, and a weight $w$, the reverse weak-type quasi-norm is
\[
\|f\|_{L^{q,\infty}(w)}
:=\inf_{\lambda>0}\lambda\,w(\{x:f(x)<\lambda\})^{1/q},
\]
where the sublevel set uses $f<\lambda$.

For $q<0<p\leq1$, we say an operator $T$ is of \emph{weak-type $(p,q)$}, if there
exists a constant $C(p,q,n)>0$ such that $\|f\|_{L^{q,\infty}}
\geq C(p,q,n)\|f\|_{L^p}$ for all $f\in L^p(\mathbb{R}^n)$.

\begin{lemma}[\cite{DouZhu2015}]
For $0<r<1$, put $r'=r/(r-1)<0$. For nonnegative measurable
functions $f,g$, the reverse H\"older inequality holds whenever the
quantities are defined:
\begin{equation*}
\int f(x)g(x)dx\ge \|f\|_{L^r}\|g\|_{L^{r'}}.
\end{equation*}
\end{lemma}

We also need the exact dual formula for a negative exponent.

\begin{lemma}\label{lem:negative-duality}
Let $q<0$ and $q'=q/(q-1)\in(0,1)$. For every positive measurable $f$,
\begin{equation*}
\|f\|_{L^q}
=\inf_{\substack{g\ge0\\ \|g\|_{L^{q'}}=1}}
\int_{\mathbb{R}^n}f(x)g(x)\,dx.
\end{equation*}
\end{lemma}
\begin{proof}
Reverse H\"older's inequality with exponents $q'$ and $q$ gives the inequality for every admissible $g$. Put $A=\int_{\mathbb{R}^{n}}f^q(x)dx$.
If $0<A<\infty$, then
$g=A^{-1/q'}f^{q-1}$ has $\|g\|_{q'}=1$ and
$$\int_{\mathbb{R}^{n}}f(x)g(x)dx=A^{1/q}=\|f\|_q.$$
If $A=\infty$, choose increasing measurable sets $E_j$ for which
$0<\int_{E_j}f^q(x)dx<\infty$ and $\int_{E_j}f^q(x)dx\to\infty$.
The same construction with $f^{q-1}\chi_{E_j}$ gives integrals
$(\int_{E_j}f^q(x)dx)^{1/q}\to0=\|f\|_q$.
\end{proof}
The strong $(p,q)$ inequality for the weighted reverse fractional integral operators can then be equivalently reformulated as a bilinear form controlled from below by a product of norms. Given $F(y)=f(y)w(y)$, define
\[
S_wF(x)=w(x)I_s(F/w)(x).
\]
Then
\[
\|I_sf\|_{L^q(w^q)}=\|S_wF\|_{L^q},
\qquad
\|f\|_{L^p(w^p)}=\|F\|_{L^p}.
\]
By Lemma \ref{lem:negative-duality}, the best strong constant equals the best constant in
\begin{equation}\label{eq:bilinear-main}
\mathcal{B}_w(F,G)
:=\iint_{\mathbb{R}^n\times\mathbb{R}^n}
|x-y|^{s-n} F(y)G(x)\frac{w(x)}{w(y)}\,dy\,dx
\ge C\|F\|_p\|G\|_{q'}.
\end{equation}
This formulation exhibits the symmetry
\begin{equation}\label{eq:bilinear-symmetry}
C_{p,q'}(w)=C_{q',p}(w^{-1}).
\end{equation}
It will be useful in the sharpness argument.

Next, we introduce the following lemma which play an important roles in the proof of the main theorems.
\begin{lemma}\label{lem:separated}
Let $\|F\|_{L^p}=\|G\|_{L^{q'}}=1$ for $q<0$. Assume $[w]_{p,q}<\infty$ and let $B=B(z,r)$. Fix $\eta\in(0,1)$. If
\begin{equation*}
\int_{B}F(y)^p\,dy\ge\eta,
\qquad
\int_{(2B)^c}G(x)^{q'}\,dx\ge\eta,
\end{equation*}
then
\begin{equation}\label{eq:sep-u-concl}
\mathcal{B}_w(F,G)
\ge c_{n,s,p,\eta}\,
[w]_{p,q}^{-1}[w^{-p'}]_{A_\infty}^{1/q}.
\end{equation}
If instead
\begin{equation*}
\int_{B}G(x)^{q'}\,dx\ge\eta,
\qquad
\int_{(2B)^c}F(y)^p\,dy\ge\eta,
\end{equation*}
then
\begin{equation*}
\mathcal{B}_w(F,G)
\ge c_{n,s,p,\eta}\,
[w]_{p,q}^{-1}[w^{q}]_{A_\infty}^{1/p'}.
\end{equation*}
\end{lemma}

\begin{proof}
We prove the first assertion. Decompose
\[
(2B)^c=\bigcup_{k\ge1}A_k,
\qquad
A_k=2^{k+1}B\setminus2^kB,
\]
and set $q_k=\int_{A_k}G(x)^{q'}\,dx$. Then $$\sum_{k\ge1}q_k\ge\eta.$$ If $y\in B$ and $x\in A_k$, then
\[
|x-y|\ge c\,2^kr.
\]
Hence
\begin{align*}
\mathcal{B}_w(F,G)
&\ge c\sum_{k\ge1}(2^kr)^{s-n}
\left(\int_B\frac{F(y)}{w(y)}\,dy\right)
\left(\int_{A_k}G(x)w(x)\,dx\right).
\end{align*}
Reverse H\"older gives
\[
\int_B\frac{F(y)}{w(y)}dy
\ge \big(\int_{B}F(y)^p\,dy\big)^{1/p}\big(\int_{B}w^{-p'}(y)dy\big)^{\frac{1}{p'}}
\ge \eta^{1/p}\big(\int_{B}w^{-p'}(y)dy\big)^{\frac{1}{p'}},
\]
and
\[
\int_{A_k}G(x)w(x)dx
\ge q_k^{\frac{1}{q'}}\big(\int_{A_k}w^{q}(y)dy\big)^{\frac{1}{q}}
\ge q_k^{\frac{1}{q'}}\big(\int_{2^{k+1}B}w^{q}(y)dy\big)^{\frac{1}{q}}.
\]
Applying the definition of the weight condition on $2^{k+1}B$ yields
\begin{align}
\mathcal{B}_w(F,G)
&\ge c\eta^{1/p}\sum_{k\ge1}q_k^{\frac{1}{q'}}(2^kr)^{s-n}\big(\int_{B}w^{-p'}(y)dy\big)^{\frac{1}{p'}}\big(\int_{2^{k+1}B}w^{q}(y)dy\big)^{\frac{1}{q}}\nonumber\\
&\ge c_{n,s,p,\eta}\,[w]_{p,q}^{-1}
\sum_{k\ge1}
\left(\frac{w^{-p'}(2^{k+1}B)}{w^{-p'}(B)}\right)^{-\frac{1}{p'}}q_k^{\frac{1}{q'}}.
\label{eq:sep-sum}
\end{align}
By the $A_\infty$ estimate \eqref{eq:Ainf-decay}, we have
\[
\frac{w^{-p'}(B)}{w^{-p'}(2^{k+1}B)}
\le C_n2^{-kn\delta_u},
\qquad
\delta_u\ge \frac{c_n}{[w^{-p'}]_{A_\infty}}.
\]
Therefore
\[
\left(\frac{w^{-p'}(2^{k+1}B)}{w^{-p'}(B)}\right)^{-\frac{1}{p'}}
\ge c_n2^{{-\frac{n\delta_u}{p'}}k},
\qquad
\text{where } 2^{{-\frac{n\delta_u}{p'}}}>1.
\]
Set $a_k=2^{-n\delta_u k/p'}$. Since $1/q+1/q'=1$ and $q<0$,
discrete reverse H\"older gives
\begin{align*}
\sum_{k\ge1}
\left(\frac{w^{-p'}(2^{k+1}B)}{w^{-p'}(B)}\right)^{-1/p'}q_k^{1/q'}
&\gtrsim_n\sum_{k\ge1}a_kq_k^{1/q'}\\
&\ge \left(\sum_{k\ge1}a_k^q\right)^{1/q}
       \left(\sum_{k\ge1}q_k\right)^{1/q'}\\
&\ge \eta^{1/q'}
2^{-n\delta_u/p'}
\left(1-2^{-n\delta_u q/p'}\right)^{-1/q}.
\end{align*}
We may choose $0<\delta_u\le1$ with
$\delta_u\gtrsim_n[w^{-p'}]_{A_\infty}^{-1}$.
Since $q/p'>0$, the last expression is bounded below by
$c_{n,s,p,\eta}\delta_u^{-1/q}
\gtrsim_{n,s,p,\eta}[w^{-p'}]_{A_\infty}^{1/q}$.
Together with \eqref{eq:sep-sum}, this proves \eqref{eq:sep-u-concl}.
The second assertion follows by exchanging $(F,p,w)$ with
$(G,q',w^{-1})$; the weight characteristic is unchanged under this exchange.
\end{proof}
The following lemma guarantees that such separated functions in Lemma \ref{lem:separated} always hold.
\begin{lemma}\label{lem:sep}
Let $0<p<1$, $0<q'<\infty$, and let $F,G\geq0$ be measurable functions on
$\mathbb{R}^n$ with
\[
\int_{\mathbb{R}^n}F(y)^p\,dy=1,
\qquad
\int_{\mathbb{R}^n}G(x)^{q'}\,dx=1.
\]
And there is a dimensional constant $\eta_n>0$
such that at least one of the following alternatives holds for some ball
$B$:
\begin{align}
&\int_B F^p(y)\,dy\geq\eta_n,
\qquad
\int_{(2B)^c}G^{q'}(x)\,dx\geq\frac12;\label{eq:alt1}\\
&\int_B G^{q'}(x)\,dx\geq\eta_n,
\qquad
\int_{(2B)^c}F^p(y)\,dy\geq\frac12.\label{eq:alt2}
\end{align}
\end{lemma}

\begin{proof}
For a nonnegative
integrable function $H$ with $\int_{\mathbb{R}^n}H(y)\,dy=1$, define
\[
c_H(r)=\sup_{z\in\mathbb{R}^n}\int_{B(z,r)}H(y)\,dy,
\qquad
r_H=\inf\{r>0:c_H(r)\geq 1/2\}.
\]
Since $H\in L^1$, absolute
continuity of the integral gives $c_H(r)\to0$ as $r\downarrow0$,
uniformly in the center. Since $\int_{\mathbb{R}^n}H(y)\,dy=1$, there holds $c_H(r)\to1$ as
$r\to\infty$. Hence $0<r_H<\infty$, and by definition every ball $B$ of
radius strictly smaller than $r_H$ satisfies $\int_{B}H(y)\,dy<1/2$.

We claim there is a ball $B_H$ of radius $2r_H$ with
$\int_{B_H}H\geq1/3$. Indeed, choose $r_j\downarrow r_H$ with
$r_j<2r_H$ and $c_H(r_j)\geq1/2$. For each $j$ there is a center
$z_j$ with $\int_{B(z_j,r_j)}H>1/3$, and
$B(z_j,2r_H)\supset B(z_j,r_j)$.

Apply this to $H=F^p$ and $H=G^{q'}$, obtaining balls $B_{F},B_{G}$ of radii
$2r_{F},2r_G$ respectively, each carrying mass $\geq1/3$. Choose $L=16$.

\emph{Case 1: $r_F\leq r_G/L$.} Then $2B_F$ has radius $4r_F<r_G$, so
$\int_{2B_F}G^{q'}<1/2$ while $\int_{B_F}F^p\geq1/3$. Thus \eqref{eq:alt1} holds
with $B=B_F$ and $\eta_n=1/3$.

\emph{Case 2: $r_G\leq r_F/L$.} Symmetrically, \eqref{eq:alt2} holds with
$B=B_G$.

\emph{Case 3: $L^{-1}r_G<r_F<Lr_G$.} Cover $B_F$, whose radius is
$2r_F$, by at most $N_n$ balls of radius $r_G/4$, where $N_n$ depends only
on $n$ and $L$. Such a bound follows by taking a maximal
$r_G/4$-separated set of centers in $B_F$: the radius-$r_G/8$ balls are
disjoint and lie in a ball of controlled radius, so their number is
bounded by a volume ratio. Since $\int_{B_F}F^p(y)dy\geq1/3$, one of these
covering balls, say $B^*$, satisfies
\[
\int_{B^*}F^p\geq\frac{1}{3N_n}.
\]
The radius of $2B^*$ is $r_G/2<r_G$, hence
$\int_{2B^*}G^{q'}<1/2$. This proves \eqref{eq:alt1} with
$\eta_n=(3N_n)^{-1}$. Interchanging $F^p$ and $G^{q'}$ proves \eqref{eq:alt2} as
well.
\end{proof}
\section{Necessity of the Muckenhoupt condition}
\begin{proposition}\label{prop:necessity}
If \eqref{eq:strong} holds with constant $C>0$, then
\begin{equation}\label{eq:necessity-quant}
[w]_{p,q}\le C_{n,s}\,C^{-1}.
\end{equation}
In particular, $w^{-p'}\in A_{1+\frac{p'}{q}}$.
\end{proposition}

\begin{proof}
Fix a ball $B=B(x_0,R)$. We first show $w^q(B)<\infty$.
For sufficiently large $N$, the set
$D_N=B\cap\{N^{-1}\le w\le N\}$ has positive measure. For $f=\chi_{D_N}$, both sides of \eqref{eq:strong} are positive and finite on $B$. If $w^q(B)=\infty$, the negative
$q$-norm of $I_sf$ is zero, contradicting (\ref{eq:strong}).
Now let $f\ge0$ be supported in $B$. For $x,y\in B$,
\[
|x-y|^{s-n}\le (2R)^{s-n},
\]
so
\[
I_sf(x)\le (2R)^{s-n}\int_Bf(y)dy.
\]
Restricting the integral to $B$ gives an upper bound for the negative norm. Hence
\[
\|I_sf\|_{L^q(w^q)}
\le (2R)^{s-n}\left(\int_Bf(y)dy\right)\left(\int_Bw(y)^{q}dy\right)^{1/q}.
\]
Combining this with the assumed lower bound yields
\[
C\le (2R)^{s-n}
\frac{\int_Bf(y)dy}{\left(\int_B f(y)^pw(y)^pdy\right)^{1/p}}
\left(\int_Bw(y)^{q}dy\right)^{1/q}.
\]
Note that the cut-off
function $f_N=w^{-p'}\mathbf1_{B\cap\{N^{-1}\le w\le N\}}$, whose integral is
always finite. If $\int_B w^{-p'}(y)dy=\infty$, the resulting ratio tends to zero, contradicting the assumed lower bound. Hence $w^{-p'}$ is locally
integrable, and using reverse H\"older's inequality, we have
\[
\inf_f\frac{\int_Bf(y)dy}{\left(\int_B f(y)^pw(y)^pdy\right)^{1/p}}
=\left(\int_Bw(y)^{-p'}dy\right)^{1/p'}.
\]
Thus
\[
\left(\int_Bw^{-p'}\right)^{-1/p'}
\left(\int_Bw^{q}\right)^{-1/q}
\le C^{-1}(2R)^{s-n}.
\]
Since $|B|^{-1/p'-1/q}=|B|^{(s-n)/n}\approx_n R^{s-n}$, division by $|B|^{(s-n)/n}$ proves \eqref{eq:necessity-quant}.
\end{proof}

\section{Strong-type bounds and sharpness}
\begin{proof}[Proof of Theorem \ref{thm:mixed}]
By homogeneity, normalize $\|F\|_p=\|G\|_{q'}=1$ and apply Lemma \ref{lem:sep}. If (\ref{eq:alt1}) holds, then
Lemma \ref{lem:separated} gives
\[
\mathcal{B}_w(F,G)
\ge c_{n,s,p}[w]_{p,q}^{-1}[w^{-p'}]_{A_\infty}^{1/q}.
\]
If (\ref{eq:alt2}) holds, then
\[
\mathcal{B}_w(F,G)
\ge c_{n,s,p}[w]_{p,q}^{-1}[w^q]_{A_\infty}^{1/p'}.
\]
Thus, we obtain
\[
\mathcal{B}_w(F,G)
\ge c_{n,s,p}[w]_{p,q}^{-1}
\min\{[w^{-p'}]_{A_\infty}^{1/q},[w^q]_{A_\infty}^{1/p'}\}.
\]
Taking the infimum over $G$ with $\|G\|_{q'}=1$ and using Lemma \ref{lem:negative-duality} proves (\ref{eq:mixed-main}).
\end{proof}

\begin{proof}[Proof of Theorem \ref{thm:characterization}]
The implication (i)$\Rightarrow$(ii) is Proposition \ref{prop:necessity}. If (ii)
holds, then \eqref{eq:Arelation} gives
$w^{-p'}\in A_{1+\frac{p'}{q}}$ and its dual weight
$w^q\in A_{1+\frac{q}{p'}}$, so both are in $A_\infty$.
The mixed lower bound of Theorem \ref{thm:mixed} is therefore strictly positive
and yields (i).
\end{proof}

In the following subsections, we first prove the lower bound, then show that the exponent cannot be improved.
\begin{proof}
Since
\[
[w^{-p'}]_{A_\infty}\lesssim_{n,p,q}[w^{-p'}]_{A_{1+\frac{p'}{q}}}
=[w]_{p,q}^{-p'}
\]
and, using \eqref{eq:dualAp},
\[
[w^q]_{A_\infty}
\lesssim_{n,p,q}[w^q]_{A_{1+\frac{q}{p'}}}
=[w]_{p,q}^{-q},
\]
Theorem \ref{thm:mixed} gives the lower bound estimate
\begin{align*}
C_*(w)
&\gtrsim
[w]_{p,q}^{-1}
\min\{[w]_{p,q}^{-p'/q},[w]_{p,q}^{-q/p'}\}\\
&=
[w]_{p,q}^{-\left(1+\max\{p'/q,q/p'\}\right)}
=[w]_{p,q}^{-\Gamma_*}.
\end{align*}
\end{proof}
Next, we show the sharpness of the exponent. Note that
\begin{equation}\label{eq:powerAp}
|x|^\delta\in A_r
\quad\Longleftrightarrow\quad
-n<\delta<n(r-1),
\end{equation}
and, as $\delta\rightarrow -n$ and $\delta\rightarrow n(r-1)$,
\begin{align}
[|x|^{-n+\varepsilon}]_{A_r}&\approx_{n,r}\varepsilon^{-1},
\label{eq:power-left-char}\\
[|x|^{n(r-1)-\varepsilon}]_{A_r}&\approx_{n,r}\varepsilon^{-(r-1)}.
\label{eq:power-right-char}
\end{align}

Let
\[
w_\varepsilon^-(x)=|x|^{(n-\varepsilon)/p'},
\qquad 0<\varepsilon<1.
\]
By \eqref{eq:power-left-char},
\begin{equation}\label{eq:Aminus}
[w_\varepsilon^-]_{p,q}
\approx\varepsilon^{1/p'}.
\end{equation}
Choose
\[
f_\varepsilon(x)=|x|^{-n+\varepsilon}\chi_{B(0,1)}(x).
\]
Since $-p/p'=1-p$,
\begin{equation}\label{eq:input-minus}
\|f_\varepsilon\|_{L^p((w_\varepsilon^-)^p)}^p
=\int_{B(0,1)}|x|^{-n+\varepsilon}\,dx
=\frac{|S^{n-1}|}{\varepsilon}.
\end{equation}
For $|x|\ge2$ and $|y|<1$, $|x-y|\le\tfrac32|x|$, hence
\[
I_sf_\varepsilon(x)
\le
\left(\frac32\right)^{s-n}|x|^{s-n}
\frac{|S^{n-1}|}{\varepsilon}.
\]
And
\begin{align*}
\int_{\mathbb{R}^n}(I_sf_\varepsilon(x))^{q}(w_\varepsilon^-(x))^{q}dx
&\ge c_{n,s}\varepsilon^{-q}
\int_{|x|>2}|x|^{(s-n)q+(n-\varepsilon)q/p'}\,dx\\
&=c_{n,s}\varepsilon^{-q}|S^{n-1}|
\int_2^\infty r^{-1-\varepsilon q/p'}\,dr\\
&\ge c_{n,s,p}\varepsilon^{-q-1}.
\end{align*}
Therefore
\begin{equation}\label{eq:output-minus}
\|I_sf_\varepsilon\|_{L^{q}((w_\varepsilon^-)^{q})}
\le C_{n,s,p}\varepsilon^{-1/q-1}.
\end{equation}
Combining \eqref{eq:input-minus} and \eqref{eq:output-minus},
\begin{equation}\label{eq:ratio-minus}
C_*(w_\varepsilon^-)
\le C_{n,s,p}\varepsilon^{-1/p'-1/q}
=C_{n,s,p}\varepsilon^{(s-n)/n}.
\end{equation}
\eqref{eq:Aminus} implies
\begin{equation}\label{eq:sharp-rho}
C_*(w_\varepsilon^-)
\lesssim [w_\varepsilon^-]_{p,q}^{-(1+p'/q)}.
\end{equation}
Hence, we obtain $C_*(w)\gtrsim [w]_{p,q}^{-\Gamma}$
with $\Gamma\ge 1+p'/q$.

On the other hand, let
\[
w_\varepsilon^+(x)=|x|^{-\frac{n}{q}+\frac{\varepsilon}{p'}},
\qquad 0<\varepsilon<1.
\]
By \eqref{eq:power-right-char},
\begin{equation}\label{eq:Aplus}
[w_\varepsilon^+]_{p,q}
\approx \varepsilon^{1/q}.
\end{equation}
The bilinear symmetry \eqref{eq:bilinear-symmetry} exchanges
\[
(p,-p',w^{-p'})\longleftrightarrow(q',-q,w^{q}),
\qquad
w\longleftrightarrow w^{-1}.
\]
For $w=w_\varepsilon^+$,
\[
\sigma=(w_\varepsilon^+)^{q}
=|x|^{-n+\varepsilon q/p'},
\]
which is just the lower-boundary family with the parameters exchanged. Applying \eqref{eq:ratio-minus} with $-p'$ and $-q$ exchanged gives
\[
C_*(w_\varepsilon^+)
\lesssim \varepsilon^{-1/p'-1/q}
=\varepsilon^{(s-n)/n}.
\]
\eqref{eq:Aplus} gives
\begin{equation*}
C_*(w_\varepsilon^+)
\lesssim [w_\varepsilon^+]_{p,q}^{-(1+q/p')}.
\end{equation*}
Hence $\Gamma\ge1+q/p'$. Together with \eqref{eq:sharp-rho},
this proves optimality of
$\Gamma_*=\max\{1+p'/q,1+q/p'\}$ and completes the proof of sharp exponent of Theorem \ref{thm:mixed}.
\section{Endpoint theorem and explicit characteristic bounds}\label{sec:endpoint}
We now consider the endpoint case $p=1$. For a positive measure $\sigma$, we use
\[
\|h\|_{L^{\frac{n}{n-s},\infty}(\sigma)}
=\inf_{\lambda>0}\lambda\,
\sigma(\{x:h(x)<\lambda\})^{\frac{n-s}{n}}.
\]
The $A_1$ characteristic is
\[
[\sigma]_{A_1}=\sup_B\frac{\frac{1}{|B|}\int_B\sigma}{essinf_B\sigma}.
\]
\begin{proof}
We first prove the sufficient. Assume $w^{\frac{n}{n-s}}\in A_1$ and put
$A=[w^{\frac{n}{n-s}}]_{A_1}$. Fix $f\ge0$ and $\lambda>0$, and set
\[
E_{\lambda}=\{x: I_{s} f(x)<\lambda\}.
\]
If \(f\) is nonzero and \(I_s f\) is finite at some point, then each superlevel set \(E_\lambda\) is bounded. Indeed, if \(\int_B f(y)dy > 0\) for some ball \(B\), then \(I_s f(x)\) grows at least like \(|x|^{s-n}\) as \(|x|\to\infty\). Thus it suffices to consider
$0<w^{\frac{n}{n-s}}(E_\lambda)<\infty$. For almost every $y$ and every $r>0$,
the $A_1$ condition gives
\[
w^{\frac{n}{n-s}}(E_\lambda\cap B(y,r))\le w^{\frac{n}{n-s}}(B(y,r))
\le Aw^{\frac{n}{n-s}}(y)v_n r^n.
\]
Set
$R_y=(\frac{w^{\frac{n}{n-s}}(E_\lambda)}{(A v_nw^{\frac{n}{n-s}}(y))})^{1/n}$.
The layer-cake identity and Tonelli's theorem give
\begin{align*}
\int_{E_\lambda}|x-y|^{s-n}w^{\frac{n}{n-s}}(x)\,dx
&=(s-n)\int_0^\infty r^{s-n-1}
  w^{\frac{n}{n-s}}(E_\lambda\setminus B(y,r))\,dr\\
&\ge(s-n)\int_0^{R_y}r^{s-n-1}
  \bigl(w^{\frac{n}{n-s}}(E_\lambda)-Aw^{\frac{n}{n-s}}(y)v_n r^n\bigr)\,dr\\
&=\frac ns(Av_n)^{\frac{n-s}{n}}w^{\frac{n}{n-s}}(E_\lambda)^{s/n}
  w^{\frac{n}{n-s}}(y)^{1/q_e}\\
&=\frac ns(Av_n)^{\frac{n-s}{n}}w^{\frac{n}{n-s}}(E_\lambda)^{s/n}w(y).
\end{align*}
Consequently, Tonelli's theorem yields
\[
\lambda w^{\frac{n}{n-s}}(E_\lambda)
\ge\int_{E_\lambda}I_sf(x)w^{\frac{n}{n-s}}(x)\,dx
\ge\frac ns(Av_n)^{\frac{n-s}{n}}w^{\frac{n}{n-s}}(E_\lambda)^{s/n}
  \int_{\mathbb{R}^n}f(y)w(y)\,dy.
\]
Division by $w^{\frac{n}{n-s}}(E_\lambda)^{s/n}$ and then taking the infimum
over $\lambda$, we obtain
\begin{equation*}
\|I_sf\|_{L^{\frac{n}{n-s},\infty}(w^{\frac{n}{n-s}})}\ge \frac ns\,v_n^{\frac{n-s}{n}}[w^{\frac{n}{n-s}}]_{A_1}^{\frac{n-s}{n}}\|f\|_{L^1(w)}.
\end{equation*}
We show necessity. First, $w^{\frac{n}{n-s}}$ is locally integrable. Indeed, for
any ball $B$, choose $N$ so that $B\cap D_N$ has positive measure,
where $D_N=\{N^{-1}\le w\le N\}$. With $f=\chi_{B\cap D_N}$,
$\|f\|_{L^1(w)}>0$ and $I_sf$ is bounded on $B$. If $w^{\frac{n}{n-s}}(B)=\infty$,
which makes the reverse weak norm zero,
contradicting \begin{equation*}
\|I_sf\|_{L^{\frac{n}{n-s},\infty}(w^{\frac{n}{n-s}})}\ge \frac ns\,v_n^{\frac{n-s}{n}}[w^{\frac{n}{n-s}}]_{A_1}^{\frac{n-s}{n}}\|f\|_{L^1(w)}.
\end{equation*}

Fix $N$ and choose a point $y\in D_N$ that is also a
Lebesgue point of $w\chi_{D_N}$ and $w^{\frac{n}{n-s}}$. Define
\[
f_\varepsilon(z)=
\frac{\chi_{B(y,\varepsilon)\cap D_N}(z)}
{\int_{B(y,\varepsilon)\cap D_N}w(x)dx},
\qquad \int_{\mathbb{R}^n} f_\varepsilon(z)w(z)\,dz=1.
\]
By Lebesgue differentiation, we obtain $\int f_\varepsilon(z)dz\to w(y)^{-1}$. For $x\in B(y,R)$,
\[
I_sf_\varepsilon(x)\le (R+\varepsilon)^{s-n}
\int_{\mathbb{R}^n} f_\varepsilon(z)\,dz.
\]
Since $\frac{n}{n-s}<0$, let $\varepsilon\rightarrow0$, this yields
\begin{equation}\label{eq:centered-A1-from-C}
\frac{1}{|B|}\int_{B(y,R)}w^{\frac{n}{n-s}}(z)dz\le v_n^{-1}(\frac ns\,v_n^{\frac{n-s}{n}}[w^{\frac{n}{n-s}}]_{A_1}^{\frac{n-s}{n}})^{\frac{n}{n-s}}w^{\frac{n}{n-s}}(y)
\end{equation}
for almost every $y$ and every $R>0$.

For an arbitrary ball $B(c,R)$ and almost every $y\in B(c,R)$,
\begin{eqnarray*}
\frac{1}{|B|}\int_{B(c,R)}w^{\frac{n}{n-s}}(z)dz
&\le& 2^n\frac{1}{|B|}\int_{B(y,2R)}w^{\frac{n}{n-s}}(z)dz\\
&\le&2^nv_n^{-1}(\frac ns\,v_n^{\frac{n-s}{n}}[w^{\frac{n}{n-s}}]_{A_1}^{\frac{n-s}{n}})^{\frac{n}{n-s}}w^{\frac{n}{n-s}}(y).
\end{eqnarray*}
Take the essential infimum in $y$, we have $w^{\frac{n}{n-s}}\in A_1$. This completes the proof of Theorem \ref{thm:endpoint-main}.
\end{proof}
\begin{remark}[The endpoint coefficient need not be sharp]
For $n=1$, $s=2$, and $w\equiv1$, then
\[\|I_sf\|_{L^{\frac{n}{n-s},\infty}(w^{\frac{n}{n-s}})}\ge \frac{1}{4}\|f\|_{L^{1}(w)}.\] The actual best constant is $1/2$.
Indeed, if $0<M=\int_\mathbb{R} f<\infty$ and $\int_{\mathbb{R}} |xf(x)|\, dx < \infty$
, Jensen's
inequality gives
\[
I_2f(x)=\int_\mathbb{R}|x-y|f(y)\,dy\ge M|x-\bar y|,
\qquad \bar y=M^{-1}\int_\mathbb{R} yf(y)\,dy.
\]
Hence $|\{I_2f<\lambda\}|\le2\lambda/M$, so
$\|I_2f\|_{L^{-1,\infty}}\ge M/2$; If $\int_{\mathbb{R}} |xf(x)|\, dx =\infty$, then
$I_2f\equiv\infty$ and the inequality is immediate. If $M=\infty$,
then $I_2f$ is infinite except possibly at one point, so the same
inequality holds in the extended sense. Conversely, for
any nonzero compactly supported $f$, the sublevel-set length is
asymptotic to $2\lambda/M$ as $\lambda\to\infty$, giving the upper
bound $1/2$.
\end{remark}

\bibliographystyle{amsalpha}

\begin{thebibliography}{10}


\vskip0.5cm
\bibitem{Buckley1993}
S.~M. Buckley,
\emph{Estimates for operator norms on weighted spaces and reverse Jensen inequalities},
Trans. Amer. Math. Soc. \textbf{340} (1993), no.~1, 253--272.

\bibitem{CarrilloEtAl2019}
J.~A. Carrillo, M.~G. Delgadino, J. Dolbeault, R.~L. Frank, and F. Hoffmann,
\emph{Reverse Hardy--Littlewood--Sobolev inequalities},
J. Math. Pures Appl. \textbf{132} (2019), 133--165.

\bibitem{ChenLiuLuTao2018}
L. Chen, Z. Liu, G. Lu, and C. Tao,
\emph{Reverse Stein--Weiss inequalities and existence of their extremal functions},
Trans. Amer. Math. Soc. \textbf{370} (2018), no.~12, 8429--8450.

\bibitem{CoifmanFefferman1974}
R.~R. Coifman and C. Fefferman,
\emph{Weighted norm inequalities for maximal functions and singular integrals},
Studia Math. \textbf{51} (1974), 241--250.


\bibitem{CruzUribeMoen2013}
D. Cruz-Uribe and K. Moen,
\emph{One and two weight norm inequalities for Riesz potentials},
Illinois J. Math. \textbf{57} (2013), no.~1, 295--323.



\bibitem{DouZhu2015}
J. Dou and M. Zhu,
\emph{Reversed Hardy--Littlewood--Sobolev inequality},
Int. Math. Res. Not. IMRN \textbf{2015} (2015), no.~19, 9696--9726.

\bibitem{DGPP2005}
O. Dragi\v{c}evi\'c, L. Grafakos, C. Pereyra and S. Petermichl,
\emph{Extrapolation and sharp norm estimates for classical operators on weighted Lebesgue spaces},
Publ. Math. \textbf{49} (2005), no.~1, 73--91.

\bibitem {FL} R. Frank and E. Lieb, \textit{A new, rearrangement-free proof of the Sharp Hardy-Littlewood-Sobolev
Inequality}. Spectral Theory, Function Spaces and Inequalities. Oper. Theory: Adv. Appl. \textbf{219} (2012), 55--67.



\bibitem {HL} G. Hardy and J. Littlewood, \textit{Some properties of fractional integral I }. Math. Zeit.  (1928), 565--606.

\bibitem{HuangTang2023}
L. Huang and H. Tang,
\emph{Optimal weak estimates for Riesz potentials},
C. R. Math. Acad. Sci. Paris \textbf{361} (2023), 1123--1131.


\bibitem{HuangTangZhang2025}
L. Huang, H. Tang and C. Zhang,
\emph{Endpoint estimates of fractional integral operators when $\alpha>n$},
Bull. Sci. Math. \textbf{199} (2025), 103569.

\bibitem{HuntMuckenhouptWheeden1973}
R.~Hunt, B.~Muckenhoupt, and R.~Wheeden,
\emph{Weighted norm inequalities for the conjugate function and Hilbert
transform},
Trans. Amer. Math. Soc. \textbf{176} (1973), 227--251.

\bibitem{Hytonen2012}
T.~P. Hyt\"onen,
\emph{The sharp weighted bound for general Calder\'on--Zygmund operators},
Ann. of Math. (2) \textbf{175} (2012), no.~3, 1473--1506.

\bibitem{HytonenPerez2013}
T. Hyt\"onen and C. P\'erez,
\emph{Sharp weighted bounds involving $A_\infty$},
Anal. PDE \textbf{6} (2013), no.~4, 777--818.

\bibitem{HytonenPerezRela2012}
T. Hyt\"onen, C. P\'erez, and E. Rela,
\emph{Sharp reverse H\"older property for $A_\infty$ weights on spaces of homogeneous type},
J. Funct. Anal. \textbf{263} (2012), no.~12, 3883--3899.


\bibitem{KassymovRuzhanskySuragan2022}
A. Kassymov, M. Ruzhansky, and D. Suragan,
\emph{Reverse Stein--Weiss, Hardy--Littlewood--Sobolev, Hardy, Sobolev and Caffarelli--Kohn--Nirenberg inequalities on homogeneous groups},
Forum Math. \textbf{34} (2022), no.~5, 1147--1158.

\bibitem{LaceyMoenPerezTorres2010}
M.~T. Lacey, K. Moen, C. P\'erez, and R.~H. Torres,
\emph{Sharp weighted bounds for fractional integral operators},
J. Funct. Anal. \textbf{259} (2010), no.~5, 1073--1097.

\bibitem{Lerner2008}
A.~Lerner,
\emph{An elementary approach to several results on the Hardy--Littlewood
maximal operator},
Proc. Amer. Math. Soc. \textbf{136} (2008), 2829--2833.

\bibitem{Lerner2013}
A.~K. Lerner,
\emph{A simple proof of the $A_2$ conjecture},
Int. Math. Res. Not. IMRN \textbf{2013} (2013), no.~14, 3159--3170.


\bibitem {L}E. Lieb, \textit{Sharp constants in the Hardy--Littlewood--Sobolev and related inequalities}. Ann. of Math. (2) \textbf{118} (1983), no.~2, 349--374.

\bibitem{Muckenhoupt1972Maximal}
B. Muckenhoupt,
\emph{Weighted norm inequalities for the Hardy maximal function},
Trans. Amer. Math. Soc. \textbf{165} (1972), 207--226.

\bibitem{MuckenhouptWheeden1974}
B. Muckenhoupt and R.~L. Wheeden,
\emph{Weighted norm inequalities for fractional integrals},
Trans. Amer. Math. Soc. \textbf{192} (1974), 261--274.
\bibitem{NgoNguyen2017}
Q.~A. Ng\^o and V.~H. Nguyen,
\emph{Sharp reversed Hardy--Littlewood--Sobolev inequality on $\mathbb R^n$},
Israel J. Math. \textbf{220} (2017), 189--223.

\bibitem{NgoNguyen2017Half}
Q.~A. Ng\^o and V.~H. Nguyen,
\emph{Sharp reversed Hardy--Littlewood--Sobolev inequality on the half space
$\mathbb R^n_+$},
Int. Math. Res. Not. IMRN \textbf{2017}, no.~20, 6187--6230.


\bibitem{Petermichl2007}
S.~Petermichl,
\emph{The sharp bound for the Hilbert transform in weighted Lebesgue spaces
in terms of the classical $A_p$ characteristic},
Amer. J. Math. \textbf{129} (2007), no.~5, 1355--1375.

\bibitem{Petermichl2008}
S.~Petermichl,
\emph{The sharp weighted bound for the Riesz transforms},
Proc. Amer. Math. Soc. \textbf{136} (2008), no.~4, 1237--1249.

\bibitem{Perez1990}
C.~P\'erez,
\emph{Two weight norm inequalities for Riesz potentials and uniform
$L^p$-weighted Sobolev inequalities},
Indiana Univ. Math. J. \textbf{39} (1990), 31--44.

\bibitem{Perez1994}
C.~P\'erez,
\emph{Two weight inequalities for potential and fractional type maximal
operators},
Indiana Univ. Math. J. \textbf{43} (1994), 663--683.

\bibitem{Sawyer1982}
E.~T. Sawyer,
\emph{A characterization of a two-weight norm inequality for maximal
operators},
Studia Math. \textbf{75} (1982), 1--11.

\bibitem{Sawyer1988}
E.~T. Sawyer,
\emph{A characterization of a two-weight norm inequality for fractional and
Poisson integrals},
Trans. Amer. Math. Soc. \textbf{308} (1988), 533--545.

\bibitem{SawyerWheeden1992}
E.~T. Sawyer and R.~L. Wheeden,
\emph{Weighted inequalities for fractional integrals on Euclidean and
homogeneous spaces},
Amer. J. Math. \textbf{114} (1992), 813--874.

\bibitem {R} M. Riesz, \textit{L'int\'egrale de Riemann--Liouville et le probl\`eme de Cauchy}. Acta Math. \textbf{81} (1949), 1--222.

\bibitem{Rubio1984}
J.~L. Rubio de Francia,
\emph{Factorization theory and $A_p$ weights},
Amer. J. Math. \textbf{106} (1984), no.~3, 533--547.

\bibitem{RegueraThiele2012}
M.~C. Reguera and C. Thiele,
\emph{The Hilbert transform does not map $L^1(Mw)$ to $L^{1,\infty}(w)$},
Math. Res. Lett. \textbf{19} (2012), no.~1, 1--7.



\bibitem {S} S. Sobolev, \textit{On a theorem in functional analysis }(in Russian). Mat. Sob. \textbf{46} (1938), 471--497.

\bibitem{SteinWeiss1958}
E.~M. Stein and G. Weiss,
\emph{Fractional integrals on $n$-dimensional Euclidean space},
J. Math. Mech. \textbf{7} (1958), no.~4, 503--514.

\bibitem{TangWang2026}
H. Tang and Y. Wang,
\emph{Sharp constants for weak estimates of the Riesz potentials},
arXiv:2608.31043v2 [math.CA] (2026).

\bibitem{Wilson2008}
M. Wilson,
\emph{Weighted Littlewood--Paley Theory and Exponential-Square Integrability},
Lecture Notes in Mathematics, vol.~1924,
Springer, Berlin, 2008.

\bibitem{Zygmund1956}
A. Zygmund,
\emph{On a theorem of Marcinkiewicz concerning interpolation of operators},
J. Math. Pures Appl. (9) \textbf{35} (1956), 223--248.

\end{thebibliography}

\end{document}